\documentclass[11pt]{amsart}
\usepackage{graphicx}
\input xypic \xyoption{curve}

\def\centerbmp#1#2#3{\vskip#2\relax\centerline{\hbox to#1{\special
  {bmp:#3 x=#1, y=#2}\hfil}}}

\newcommand{\C}[1]{{\mathcal#1}} 

\theoremstyle{plain}

\newtheorem{conj}{Conjecture}[section]

\theoremstyle{definition}
\newtheorem{defin}{Definition}[section]

\theoremstyle{definition}

\theoremstyle{remark}
\newtheorem{rem}{Remark}[section]

\usepackage{fancyhdr}
\begin{document}

\title{On Prime Numbers and The Riemann Zeros}         
\author{Lucian M. Ionescu}
\address{Department of Mathematics, Illinois State University, IL 61790-4520}
\email{lmiones@ilstu.edu}

\begin{abstract}
The current research regarding the Riemann zeros suggests the existence of a 
non-trivial algebraic/analytic structure on the set of Riemann zeros \cite{Marco, Entropy, Ford}.
The duality between primes and Riemann zeta function zeros suggests some
new goals and aspects to be studied: {\em adelic duality} and the {\em POSet of prime numbers}.

The article presents computational evidence of the structure of the imaginary
parts $t$ of the non-trivial zeros of the Riemann zeta function $\rho=1/2+it$,
called in this article the {\em Riemann Spectrum},
using the study of their distribution, as in \cite{Ford}.

The novelty represents in  considering the associated characters $p^{it}$,
towards an algebraic point of view, than rather in the sense of Analytic Number Theory.
This structure is tentatively interpreted  in terms of adelic characters,
and the duality of the rationals.

Second, the POSet structure of prime numbers studied in \cite{LI-Primes},
is tentatively mirrored via duality in the Riemann spectrum. 
A direct study of the convergence of their Fourier series, along Pratt trees, is proposed.

Further considerations, relating the Riemann Spectrum, adelic characters and distributions, 
in terms of Hecke idelic characters, local zeta integrals (Mellin transform) and 
$\omega$-eigen-distributions, are explored following \cite{Kudla}.
\end{abstract}

\maketitle

%

\section{Introduction}       
One of the most famous and important problem in mathematics is the Riemann Hypothesis (RH),
not only because of the multitude of conditional theorems depending on its validity,
but also because of the connections with multiplicative number theory and
related areas of mathematics.

Deeply connected to mathematical-physics, e.g. with 
Riemann gas, quantum chaos, random matrices 
etc. (see also \cite{LI:RemPNT}), 
the zeta function zeroes are ``dual'' with the prime numbers \cite{Mazur-Primes,Marco}, 
at least in the sense that the Riemann-Mangoldt exact formula 
resembles a Poisson summation formula.
Understanding this duality would clarify the ``origin'' of these zeroes, 
beyond our capability to compute them and count them precisely,
perhaps shedding light on the ``right approach'' to prove the RH.


In a different direction, 
a natural partial order on the set of prime numbers (POSET) was recently (re)discovered \cite{LI-Primes},
being related to the concept of Pratt tree \cite{Pratt} (1975) used in primality tests.
It is therefore natural to ponder on a possible mirroring phenomenon,
as a yet to be determined structure on the set of Riemann zeros.
The present article is a preliminary investigation in this direction,
exploring the adelic duality as a theoretical avenue supporting 
the existence of such a structure, and its relation to the POSet of primes.

\vspace{.2in}
In what follows, the imaginary parts of the non-trivial zeros of the Riemann zeta function
are called the {\em Riemann spectrum} \cite{Mazur-Primes} or {\em zeta eigenvalues}
(anticipating a relation with the adelic $\omega$-eigendistributions).

On the side of the prime numbers, 
the author introduced a partial order derived conceptually from the internal symmetry structure of basic finite fields \cite{LI-Primes}:
$q<<p$ if $q|p-1$, i.e. $Aut_{Ab}(Z/pZ, +)$ has a $q$-order symmetry (generator).
This leads to a hierarchy on primes, with a rooted tree associated to each prime.

For example
$$t(181)=B^+(\bullet, t(2), t(2), t(3), t(3), t(5)), \ t(3)=B^+(\bullet), 
\ t(5)=B^+(\bullet, \bullet).$$

$$\diagram
t(181)= & & & \circ \dllto \dlto \dto \drto \drrto & \\
\phi(181)=2^2\cdot 3^2\cdot 5 
  & \overset{2}{\bullet} & \overset{2}{\bullet} & \overset{3}{\bullet} \dto & 
                    \overset{3}{\bullet} \dto & \overset{5}{\bullet} \dto \drto \\
\phi(3)=2, \ \phi(5)=2^2 & & & \overset{2}{\bullet} & \overset{2}{\bullet}
 & \overset{2}{\bullet} & \overset{2}{\bullet}
\enddiagram$$
where $B^+$ is the operator which adjoins a common root to the input rooted trees.
Of course the labels ($2,3$ and $5$) are not necessary; the corresponding
``prime values'' of the labels can be reconstructed recursively, from terminals going up.

In this way the POSet of primes is a connected directed graph, with no 
directed cycles, defining a ``flow'' with sink at $2$. 

On the other side of the primes-zeros duality, 
some progress was achieved in the understanding of the 
imaginary part of Riemann zeta function zeros (ZZ)
$\rho_n=1/2+2\pi i \gamma_n$, assuming of course the Riemann Hypothesis is true.
Namely, in \cite{Ford, Ford2} the set of imaginary parts of 
(non-trivial) zeta zeroes $\gamma=\{\gamma_n\}_n$ is investigated by looking at the changes
under multiplication by real numbers $L_\alpha(\gamma_n)=\{\alpha\gamma_n\}$, 
after reduction to the unit interval ($\{\ \}$ denotes the fractional part).

The work of \cite{Ford}, combined with the idea of duality and
the additional POSet structure on primes,
suggests the existence of additional structures on them.
More precisely, since the duality is of adelic origin, 
with a Fourier transform correspondence between
primes and zeros, a similar POSet structure is expected to exist on them.

We will present next, in more detail, the current known results regarding the 
distribution of the Riemann spectrum.

Our reinterpretation of their direct study of the imaginary parts, in terms of characters $p^{it}$,
opens an avenue for investigation using Gauss and Jacobi sums, with a possibility
to relate them to the Weil zeros, in a follow-up paper.
In this present study we will focus on exploring the adelic interpretation
and the relation with the POSet structure of the prime numbers.

\section{The Distribution of the Riemann Spectrum}\label{S:RiemannSpectrum}
It becomes more and more evident that not the zeros $\rho_n$ 
are the fundamental,
but their imaginary parts $\gamma_n$, 
with $n^{-1/2}$ a normalization factor, analog to the
one present in the discrete Fourier transform for finite Abelian groups $Z/nZ$.

It is known that the fractional parts of $\{\gamma_n\}$ are uniformly distributed on 
the unit interval \cite{Ford2}, p.1.
This result can be interpreted as dual to Dirichlet Theorem regarding the uniform
distribution of prime numbers on the ``finite circles'' $Z/pZ$.
Here the result refers to the Riemann zeros instead,  
and the finite cycle $Z/pZ$ is replaced 
by the normalization of the unit circle $R/Z \cong T=S^1$.

This approximation result is formulated as the integration pairing 
with a function of an appropriate class \cite{Ford2}, p.2, Eq. (1.1):
\begin{equation}\label{E:Ford}
\sum_{0<\gamma\le T} f(\alpha\gamma)=N(T) \int_T f(x)dx + o(N(T),
\end{equation}
\begin{equation}\label{E:NT} 
N(T)=\frac{T}{2\pi}(\log \frac{T}{2\pi}-1)+O(log T).
\end{equation}
Ford and Zaharescu found a lower order term of the approximation present
{\em only} when $\alpha=a/q\cdot \log p /2\pi$ with $r=a/q$ rational
$$  \sum_{0<\gamma\le T} f(\alpha\gamma)=N(T) \int_T f(x)dx
+ T\int_T f(x) g_\alpha(x)dx + o(T)$$
where:
\begin{equation}\label{E:distribution}
g_\alpha(x)=-\frac{\log p}{\pi} Re[ \sum_1^\infty \frac{e^{-2\pi i qkx}}{p^{ak/2}}].
\end{equation}
It is important to note that re-shuffling the imaginary parts by an irrational 
scaling factor $\alpha\ni R$ keeps their distribution on the unit interval uniform
(and $g_\alpha\cong0$).

We interpret their result as a clear indication of a hidden structure on the 
imaginary parts of zeta zeros.

\begin{rem}
We also compare their approximation formula for {\em counting zeros}
with the well-known Riemann-Mangoldt exact formula for the
{\em counting of primes}, with emphasis on the sum over Riemann zeros
$$\sum_{\rho} x^\rho/\rho +\sum_1^\infty x^{-2k}/(-2k)
=x-\sum_{n\le x} \Lambda(n) - \zeta'(0)/\zeta(0).$$
as in the old formula of Landau (\cite{Ford}, p.1):
$$\sum_{0<\gamma\le T} x^\rho=-\frac{T}{2\pi} \Lambda(x)+O(T).$$
The presupposed symmetry between primes and zeros, warranted by this yet unknown 
{\em explicit} duality, should reflect in an analog of the RM-exact formula,
yielding a count of zeros, with an appropriate measure,
in terms of the prime numbers.

Note that characteristic functions of intervals in $T$ are among the admissible functions $f$
\cite{Ford2}, p.2., 
in the above distributional approximation formula (weak convergence of measures).
\end{rem}
Our goal is to find a hypothetical ``direct sum'' decomposition on the set of zeta zeros, as 
eigenvalues belonging to different ``sectors'' of the regular {\em adelic} representation.
This would be a sort of ``localization'', a reflection of a {\em local-to-global principle}
applied to the zeros themselves.

To get a feel of the above ``resonance'' effect (role of an eigenvalue), 
of shuffling the zeros via scaling, 
we reproduce the distributions from \cite{Ford}, Fig.1, p.4.

For $p=2$ and therefore $\alpha=\log 2/2\pi$, the distribution of $\gamma_n$,
i.e. of the imaginary parts of the Riemann zeros modulo $2\pi$ and rescaled to the unit interval, 
has a pick at $0.5$.
\begin{figure} [h!]
\includegraphics[width=4in,height=2in]{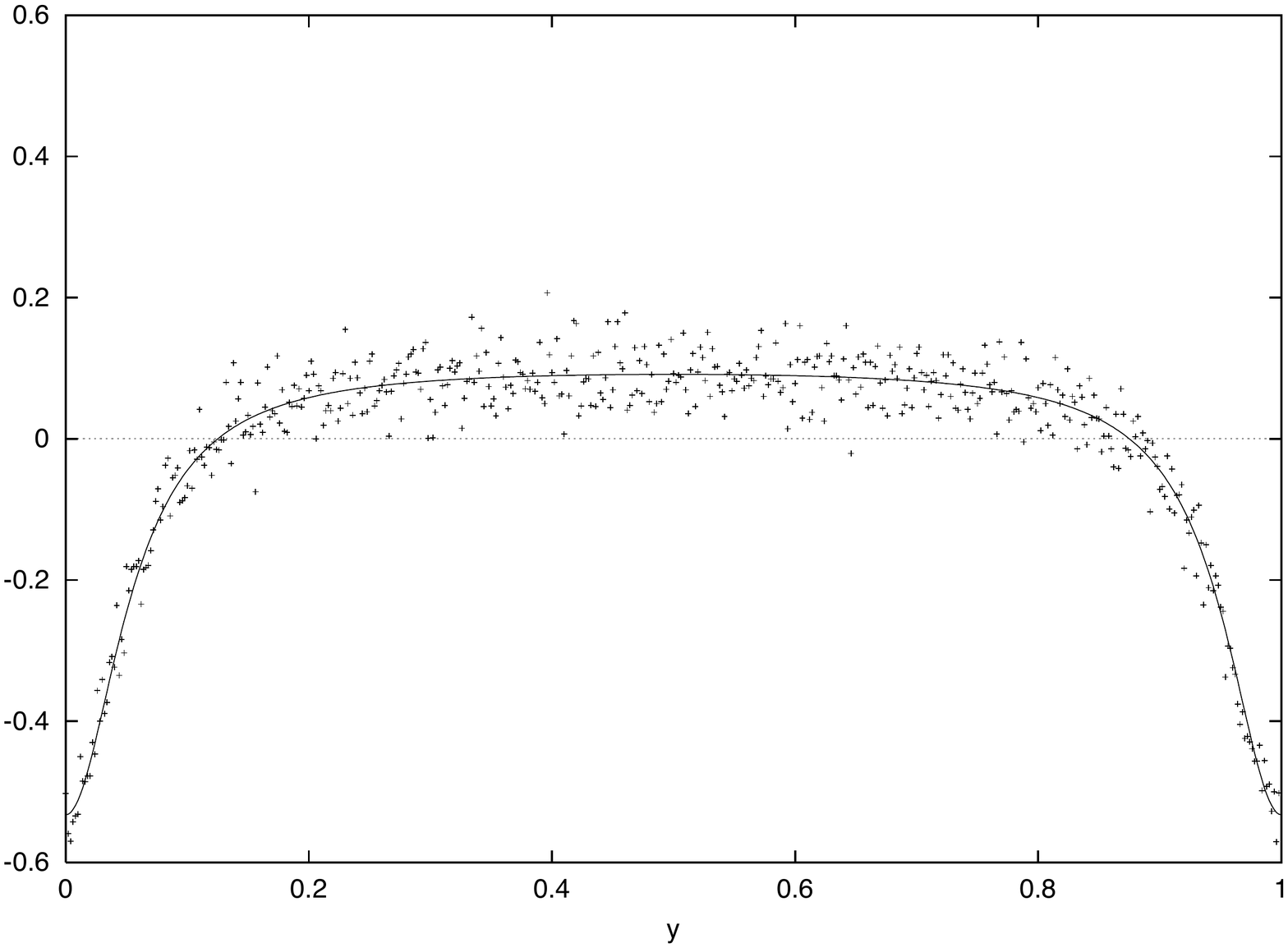} 
\vspace{-.2in} \caption{$p=2, \ q=1$}
\end{figure}
Similar non-trivial distributions occur for $p=5$ and $q=3$ (interpreted as a compression factor),
shown in Fig.2.
\begin{figure}[h!]
\includegraphics[width=4in,height=2in]{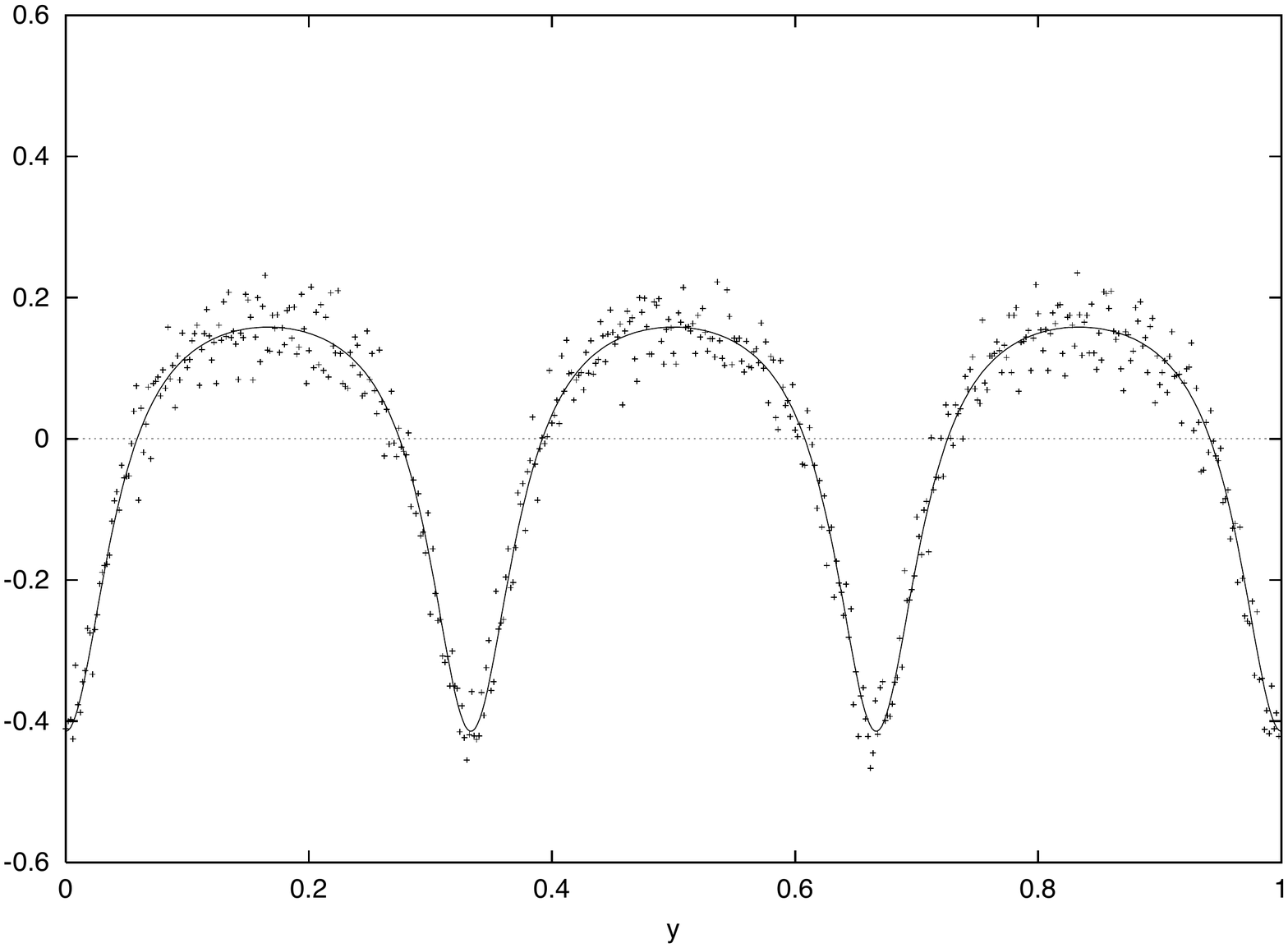} 
\vspace{-.2in} \caption{$p=5, \ q=3$}
\label{Ford-diagram-p=2}
\end{figure}

On the other hand when mixed using a composite number $p=6$, 
the uniform distribution is obtained in Fig.3.
\begin{figure} [h!]
\includegraphics[width=4in,height=2in]{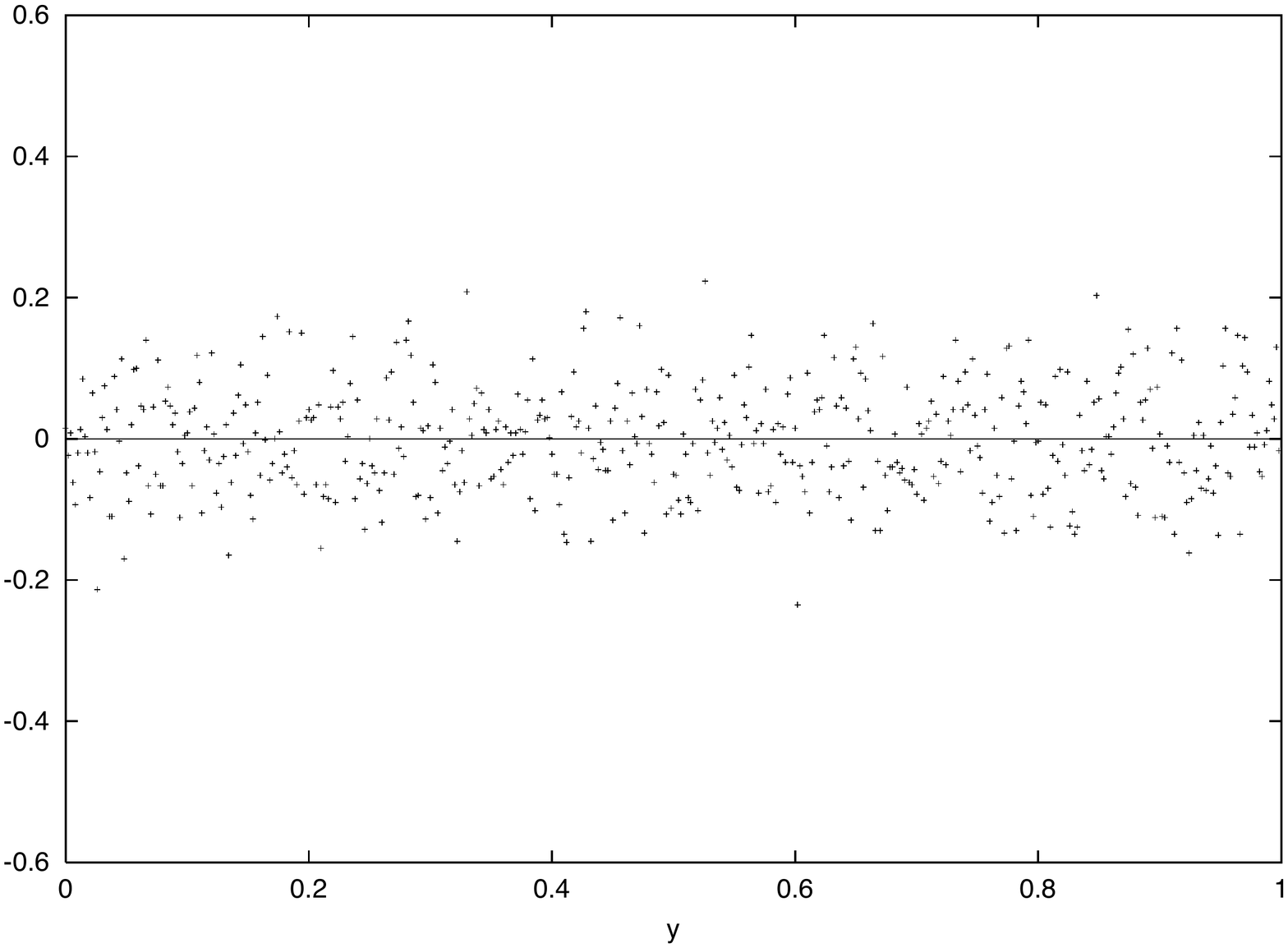} 
\vspace{-.2in} \caption{$p=6, \ q=1$}
\end{figure}

Before explaining further in \S \ref{S:p-adicSectorsOfRSpec}
the significance of these diagrams, 
we will briefly review the some facts regarding adelic duality
and Fourier transform at the level of distributions.

\section{On Quasi-crystals and Fourier Transform}\label{S:QQFTD}
Fourier transform relating the $\delta$-basis of primes in $Q$ with a dual basis of
characters of $\hat{Q}$, is reminiscent of the old hunch of crystal basis \cite{Dyson-FB},
also suggested by the Odlyzko's computation \cite{Odlyzko},
as reported by Baez \cite{Baez-qq}.

\subsection{Quasicrystals}\label{S:qq}
The fact observed by \cite{Odlyzko} (1989) and 
rediscovered my Matt McIrvin around 1990s, is that summing 
exponentials $exp(i\gamma_n x)$ yields an approximation of the
distribution $\C{D}_P=\sum_{p\in P, n\in N} \log p\ \delta_{\log p^n}$.

\begin{figure} [h!]
\includegraphics[width=4in,height=4in]{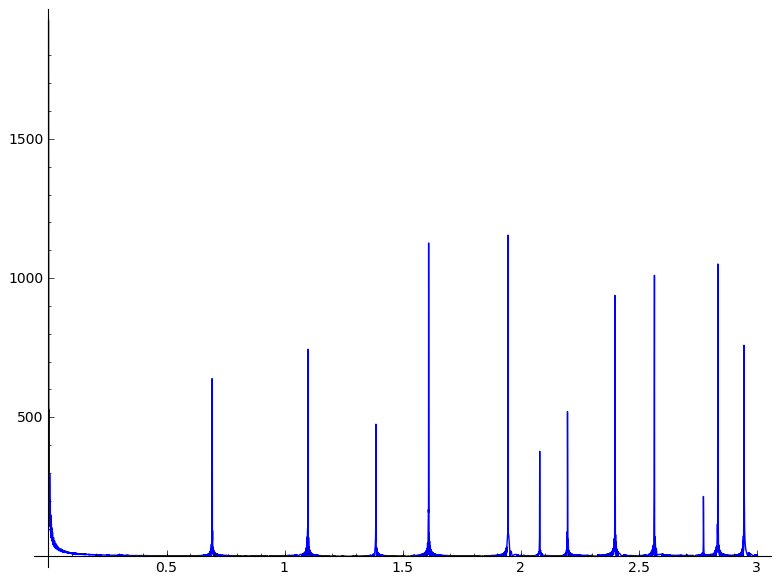}
\caption{From Baez \cite{Baez-qq}, p.4}
\label{Baez-diagram}
\end{figure}

In the following {\em probabilistic framework}:
$$\diagram
(g_P, d\mu) \rto^{Exp}_{iso} \drto_{\int} & (N^\times,\cdot, dn) \dto^{\ln(n)}
& d\mu=\log(p)dp\\ 
& (R,+, dx) & Exp(\sum k_pX_p)=\prod_p p^{k_p}.\\
\enddiagram$$
with the ``exponential'' an isomorphism of groups which
allows to interpret the Chebyshev's function as an integral
($\psi(n)=\int_P k_n(p) d\mu(p)$),
the basis of primes in $g_P=ZP$ (the corresponding $Z$-module)
has a natural weight (measure) $\log p$, so it is natural to consider
the measure $\C{M}_P=\sum_p \log p\ \delta_p$ as the counting measure
corresponding to $\C{D}_P$.

\subsection{Fourier Transform of Distributions}\label{S:FTD}
The duality between primes and zeta zeros is presently known at the level
of {\em Fourier transform of distributions},
of the corresponding characteristic sets \cite{Mazur-Primes}, Part IV.

With our adapted notation (loc. cit.: $T\leftrightarrow H, P=\Psi$), we have:
$$\diagram 
T(t)=\sum\limits_{n=0}^\infty \cos(\theta_n t) \rrto^{\quad w-lim} & &
\C{D}_P(t)=\iint\limits_{P\times N} \frac1{\sqrt{p}^n}d_{\log p^n} \dllto_{FT} 
\dto^{Duality}\\
P(\theta)=\mathop{\sum\sum}\limits_{(p,n)\in P\times N}
\frac{\log p}{\sqrt{p}^n} \cos(n\log p \ \theta) \rrto^{\quad w-lim} & & 
\C{D}_\Theta(\theta)=\sum\limits_{n=0}^\infty d_{\theta_n}(\theta) \ullto \uto
\enddiagram$$
where $d_x=1/2(\delta_{-x}+\delta_x)$ is the symmetrized Dirac delta distribution,
so that $\hat{d_n}=\cos(x)$, 
and the $(P,d\mu)$ measure space includes the weight $\log p$ in the measure.

The diagonal arrows represent taking the Fourier transform of distributions,
while the weak limits $w-\lim$ correspond to comparing the 
``peaks'' of finite approximations with the appropriate delta functions, 
as explained below.

\subsubsection{From Primes to Zeros}
The following diagram represents an approximation of the 
Fourier transform $\hat{\C{D}_P^C}=
\sum_{p,n, p^n<C} \frac{\log p}{\sqrt{p}^n} d_{\log p^n}(t)$, 
using a cut-off $C=3$ (similarly $P_C(\theta)$ is the partial sum 
with terms $p^n\le C$).
%
\begin{figure} [h!]
\includegraphics[width=6in,height=4in]{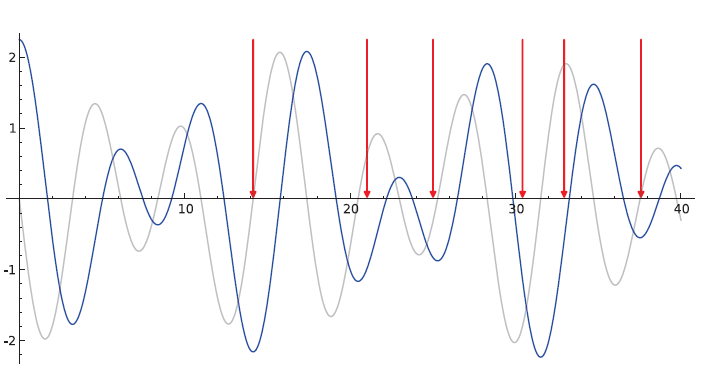}
\caption{From \cite{Mazur-Primes}, part IV, p.107}
\label{Mazur-diagram1}
\end{figure}
The peaks of the Fourier Transform seem to converge to the zeta zeros $\theta_n$
(see Fig.2).
In the distributional sense:
$$\lim_{C\to \infty}\hat{\C{D}_P^C}(\theta)= \lim_{C\to \infty} \sum_{p^n<C} 
\frac{\log p}{\sqrt{p}^n} \cos(n\log p \ \theta)
=\sum_n d_{\theta_n}(\theta)=\C{D}_{\Theta}(\theta).$$

\subsubsection{From Zeros to Primes}
In the other direction, from the Riemann spectrum to the primes,
the Inverse Fourier Transform yields in a similar way, via approximations (see Fig.3):
$$\hat{\C{D}_{\Theta}}(\log s)=\sum_n \cos(\theta_n \log s)\overset{w}{=}
\sum_{p\in P, n\in N} \log p \ d_{p^n}(s)=
\iint\limits_{P\times N} d_{p^n}(s).$$
%
\begin{figure} [h!]
\includegraphics[width=6in,height=4in]{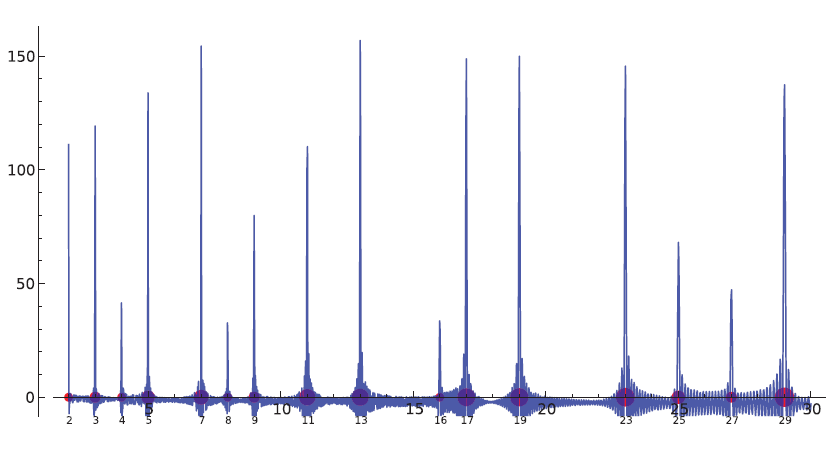}
\caption{From \cite{Mazur-Primes}, part IV, p.117}
\label{Mazur-diagram2}
\end{figure}
The maxima of the $C=1000$ approximation $T_C(\theta$ 
are localized at power of primes, the limit being $\C{D}_P\circ e^t$.
With $s=e^t$ we obtain $\hat{\C{D}_\Theta}=\C{D}_P$, as expected.

\section{The Adelic and Idelic Characters}
Since Tate's thesis \cite{Kudla}, the modern tool to understand analytic continuation and 
functional equations of L-functions proved to be harmonic analysis on the adeles;
to understand their instrumental role,
a brief review of adelic (additive) and idelic (multiplicative) characters is in order.

It is well known that the finite primes analog of Pontryagin duality 
$Z\to R\to T$ is $Q\to A\to \hat{Q}$, 
where $Q=(Q,+)$ is viewed as an additive group, $A=A_Q$ denotes the adeles of the 
field $Q$, and $\hat{Q}$ is the dual group of additive characters
\cite{Conrad} (\cite{Gelfand}, Ch.3).

The additive characters decompose as a product over the places of $Q$ (its completions):
$$\psi(x)=\prod \psi_p(x), \ \psi_\infty(x,y)=e^{2\pi i <x,y>_\infty}, x,y\in R\ 
\psi_p(a,b)=e^{<a,b>_p}, a,b\in A,$$
$$<,>:A\times A^*\to T, \ <a,b>=\sum_{q\in \bar{P}} \{a_q b^*_q\}_p,\ a,b\in A$$ 
where $<,>$ is the ``Hermitean inner product'' on adeles, 
``gluing'' $A$ and its Hermitean dual $A^*$, in analogy with the Hilbert space
formalism\footnote{Reinterpreting a Hermitean inner product as a bilinear mapping,
by conjugating scalar multiplication on the dual Hilber space.},
$\bar{P}=P\cup \{\infty\}$ is the set of places of $Q$, finite in infinite,
while $\{a_p\}_p:Q_p\to Q_p/Z_p$ is the principal part of the p-adic component
$a_p$ of the adele $a$.
Since there are only finitely many components with non-trivial principal parts,
the sum if finite. 
It is viewed as an element of the {\em adelic infinite torus}
$$T=R/Z \oplus \prod_{p\in P} Q_p/Z_p\cong S^1\oplus Q/Z.$$
Note that the orientation on the {\em rational circle} $Q/Z$
corresponds with the opposite orientation on the {\em real circle} $R/Z$,
under the Hermitean conjugation on the adeles.
It ensures the kernel of the character map $a\to \psi_a(b)$ is precisely 
the diagonal rationals.
It also corresponds to Artin's product formula of norms, archimedian and 
non-archimedian (product being 1).

\begin{rem}
The additive characters of the rational adeles are given by the following composition
$$\psi(a,b)=exp\circ \{\ \} \circ m_A:A\times A \to A \to T \to S^1,$$ 
where $\{\ \}$ is the projection onto the principal part and $m_A$
is the componentwise multiplication.
They form an affine space with structure space 
the multiplicative ``translations'' $L_a(b)=ab$:
$$\psi_a\circ L_b=\psi_{ab}.$$
Their restrictions to $Q$ yield the additive characters of $Q$.
\end{rem}

The multiplicative characters of $Q^\times$ are the obvious ones, from its structure
as a free abelian group with basis $P$.
The adelic multiplicative characters $\omega:A^\times/Q^\times$,
trivial on $Q^\times$, are called Hecke characters.
They have the form $\omega(x)=\chi(x) |x|^s$, with $\chi$ a Dirichlet character. 
The finite order Hecke characters are the Dirichlet characters.

With the adeles as the natural framework for a prime-zeros duality,
a natural question arises:
``Is there a p-adic decomposition of the Riemann Spectrum, 
a ``homological'' analog of the Euler products of the L-functions?''

\section{In Search of p-adic Sectors of the Riemann Spectrum}\label{S:p-adicSectorsOfRSpec}
Returning to the findings from \cite{Ford}, consider the reals acting on the imaginary part
of the zeros $\gamma$ via dilations $L_\alpha(x)=\alpha \gamma$, and then plotting the 
distribution of their fractional (principal) part.

When $\alpha$ is real the distribution is unchanged: uniform.

If $\alpha=a \log p/2\pi q$ ($a/q$ rational) 
the distribution $g_\alpha(x)$ (Equation \ref{E:distribution})
has global minima at the points $k/q, k=0,...,q-1$, corresponding to the 
$q$-roots of unity (\cite{Ford}, p.3):
$$g_\alpha(k/q)=\frac{\log p}{\pi (1-\sqrt{p}^a)}<0, \ k=0,...,q-1.$$
This means ``there is a shortage of zeros '' with imaginary parts having the 
fractional part close to $k/q$.

\subsection{Relation to Adeles: Start-up Conjectures}
This suggests a natural working hypothesis. But first, some handy terminology:
\begin{defin}
The imaginary parts of the Riemann zeta zeros are called the {\em zeta eigenvalues},
and denoted $\lambda_n=\gamma_n/(2\pi)$.
\end{defin}
\begin{conj}
The zeta eigenvalues are ``adelic in nature'' corresponding to prime numbers 
({\em adelic sectors}), 
while the trivial zeros correspond to the infinite prime Euler factor
\footnote{When choosing the Gaussian distribution
as a test function, it involves the gamma function \cite{Garrett-REF}.}
\end{conj}
The second part is compatible with the known structure of the zeta function:
$$\zeta=\xi\cdot 1/\Gamma$$
with $\xi$ the {\em complete zeta function}, 
corresponding to the non-trivial zeros \cite{Garrett-REF}.
\begin{rem}
This decomposition should be interpreted in the context of zeta integrals
of the rational adeles, which have a direct product decomposition
in terms of local zeta integrals \cite{Kudla}.
\end{rem}
In order to better understand this formula, we recast it in a slightly
different form.

First recall that $A\cong Q + [0,1) \times \prod Z_p$.
Also the correspondence between $[0,1)$ ($R/Z$) and $S^1$ ($p$-roots of unity
$Q_p/Z_p$; recall $Q/Z\cong \prod Q_p/Z_p$) requires a normalization by $2\pi$.

So, introduce the variable $\lambda=\gamma/(2\pi)$, 
suggested by the zeros count formula Equation \ref{E:NT}.
\begin{rem}
We take a moment to comment on its interpretation of a quantity of information.

$$N(2\pi x)=x\log x-x +O(\log x)=\int_0^x \log t dt + O(\log x)$$
is reminiscent of a quantity of information of the distribution of zeta eigenvalues
In the case of a finite probability distribution $Z=\{p_n\}$ 
($\sum p_n=1$ are probabilities, not prime numbers):
$$I(x)=\sum_{p_n<x} \log p_n, H(x)=\sum_{p_n<x} p_n \log p_n 
\overset{x\to \infty}{\to} H(Z) \ Shannon \ entropy.$$
This is compatible with the interpretation of primes as energy levels 
in the Riemann gas model, {\em except it refers to zeta eigenvalues}!

This formula counting the number of zeta zeros should have an ``exact form''
analog to the Riemann-Mangoldt explicit / exact formula for its ``dual'' 
counting function $\pi(x)$.
This suggests to look at a similar ``recounting measure'', 
analog of $\psi(x)$, but which requires to unravel the p-adic powers
grading of the zeta eigenvalues. 
\end{rem}
The ``multiplicity'' $a\log p=\log (p^a)$ corresponds to the 
$a^{th}$ degree of a p-adic number (Frobenius shift in $Z_p$).

We conjecture that $q$ must be relatively prime to $p$ (different), 
or otherwise the distribution is again uniform.

\begin{conj}
The distribution of {\em zeta eigenvalues} is unchanged under dilations by 
$2\pi \alpha=\log n/q$ when $n$ is composite (e.g. diagram with $n=6$);
\end{conj}
The other ``trivial'' case is conjectured to be when $p$ divides $q$,
since we believe that the zeta eigenvalues in the $p$-adelic sector 
are related to its roots of unity, which correspond to divisors of
$p-1=|Aut(F_p,+)|$ \cite{p-adic} \footnote{The symmetries of the tangent space 
to the abelian quantum group $Z_p$ \cite{LI-presentations}.}.

\begin{conj}
The zeta eigenvalues ``eliminated'' by the dilation operator corresponding
to $2\pi \alpha /= \log p/q$, correspond to the roots of unity of $Z_p$,
i.e. the divisors of $p-1$.
\end{conj}
This is supported by the number of minima and maxima 
of the resulting distribution: $q|p-1$,
which are exactly the descendents of $p$ in the POSet of primes.

Before venturing into a recursive chase of roots of unity corresponding to the
Pratt tree of $p$ (symmetries of symmetries of symmetries ... etc. \cite{LI-presentations})
we should obtain more information on the above ``leads''.

\subsection{Additional examples and diagrams}
Regarding the parameters $a$ and $q$ of the scaling operator with $2\pi \alpha=a/q\cdot \log p$
it seems that the fundamental case is that mentioned by Rademacher \cite{Rademacher},
when $a=1$ and $q=1$.
A different $q$, like $q=3$ in Fig.1, has the effect to wrap the distribution $q$-times.

The real problem is, if there is a way to separate the zeros in families associated
to various primes, in a similar way the rational circle $Q/Z$ decomposes in the adelic duality.

\section{Primes versus Zeros: a Theoretical Exploration}\label{S:distributions}
The above ideas suggest a few {\em questions to keep in mind} : 
1) Do the zeta eigenvalues correspond to Hecke characters? (Teichmuler!?)
2) Are they local integrals of Dirichlet characters? (Hecke finite order; p-adic?).

And what about a p-adic version of the distribution function $g_\alpha(t)$,
summing p-adic characters $exp(2\pi i \{ qt\}_p)$;
does it localizes at some of the zeroes themselves?

Indeed, for $2\pi \alpha=\log p$ the distribution can be rewritten as:
$$g_\alpha(t)=\frac{\log p}{\pi} \frac{1}{1-p^{1/2}e^{2\pi i t}}$$
and seems to {\em interpolate the zero eigenvalues} (their projection on
the principal part).
Modulo the normalization coefficient $\log/\pi$, it resembles the 
value at $p$ of $L_p(s,\chi)$, 
an L-function associated to a Hecke character $\chi$ \cite{Kudla}, p.5.
($v=p, q_v=p$).
Recall that any Hecke character has the form $\omega_s=\omega\cdot |\ |^s$,
with $\omega(x)=\chi(u), x=atu$ of finite order determined by a Dirichlet character
\cite{Kudla}, p.3.
Then $L(s,\omega)=L(0,\omega\omega_s)$ (loc. cit. p.6).
If $t_v(\omega)=p$ then for $s=\rho$ (index of the zero omitted):
$$L_p(\rho, \omega)=L_p(0,\omega_\rho)=(1-t_v(\omega_\rho))q_v^{-\rho})^{-1}=
\frac{1}{1-p^{1/2}e^{2\pi i \gamma\log p/(2\pi))}}
=\frac{1}{1-\sqrt{p} e^{2\pi i \lambda}}.$$
Now $\omega_\rho$ is unramified and only depends on the coset
$s+2\pi i / \log p Z$ (see Eigendistributions \cite{Kudla}, p.6).
This seems to be compatible with the form of the distribution $g_\alpha(t)$
of $\log p \cdot \gamma_n /(2\pi)$
i.e. suggesting $i\gamma_n \in 2\pi / \log p \in Z$ (or ``close'') for a subfamily
of zero eigenvalues (some values $n$):
$$s=\rho, \quad [s]= \frac12 + i\gamma +\frac{2\pi i}{\log p}Z=
\frac12 +\frac{2\pi i}{\log p} (\gamma\frac{\log p}{2 \pi}+Z).$$

A study in more depth the theory of $\omega$-eigendistributions,
in view of the relation with local zeta integrals, is mandatory; 
for example compare the Mellin transform (local zeta integral) from 
\cite{Garrett-REF}, P 2.2.1, p.5 (author's notation: $\C{Z}$ is the 
periodicization operator, MT is Mellin transform and $f$ is the Gaussian distribution
as a canonical choice of Schwartz function) with formula (3.6) from 
\cite{Kudla} p.10, and the GCD interpretation (Remark, p.11):
$$\zeta(s)=\frac{MT(\C{Z}(GD))}{MT(f)} \quad \leftrightarrow \quad
L(s,\omega)=\frac{z(s,\omega)}{z_0(s,\omega)}=\frac{z(s,\omega;f)}{z_0(s,\omega),f>}.$$
In fact, the local zeta integral for $p=\infty$ is Mellin transform:
$$z(s,\omega_s;f)=\int_{R_+} f(x) \omega_s(x) d\mu(x), \ \omega_s(x)=|x|^s, 
d\mu(x)=dx/x\ (Haar \ measure),$$
and the archimedean case example from \cite{Kudla}, p.13,
seems to the author to just prove the above statement (Proposition 3.5).
It also points to another important property of the Gaussian distribution:
it is a basis vector (eigendistribution theory).

Since global zeta integrals have product form representations \cite{Kudla}, p.17,
we speculate that the zeta function has a product decomposition
in terms of the L-functions corresponding to its local zeta integrals:
\begin{equation}\label{E:zeta-L}
\zeta(s)=L_\infty(s,\omega_\infty) \prod_p L_p(s, \omega_p).
\end{equation}
Here we can dispense of the need for analytic continuation either by
considering the {\em alternating zeta function} \cite{Wiki:RZF}:
$$(1-2^{1-s})\zeta(s)=\sum_{n=1}^\infty (-1)^{n-1} n^{-s},$$
or by interpreting the equality in the sense of distributions (zeta integrals).
The first factor is proportional to the $\Gamma$ function,
and responsible for the trivial zeros.
Recall that 
$$L_p(\rho_n,\omega_p)=L_p(0, \omega_{\rho_n}\omega), 
\quad \omega_{\rho_n}\omega=\chi_{N(n)} |\ |^{\rho_n}.$$
On the other hand (primes), the Riemann exact formula
$$\psi(x)=x-\sum_\rho x^\rho/\rho -\sum_{-2n} x^{-2n}/(-2n)-\log (2\pi)$$
corresponding to the logarithmic derivative of the zeta function,
leads to the distribution equation:
$$\psi'(x)=1-\sum_\rho x^\rho - \sum trivial \ zeros,$$
while the LHS is a sum of Dirac's delta distributions:
$$\psi'(x)=1-\sum_{p\in P, k\in N} \log p \ \delta_{p^k}.$$
By analogy, the similar formula counting zeta zeros like in \cite{Ford},
interpreted in the sense of distributions is formally similar to:
$$\sum_\rho c_\rho \delta_\rho=\sum_p g_p(t) <-> \sum_p \log p L_p(\rho,\omega),$$
since $g_p(t)$ seem to interpolate the scaled zeta eigenvalues.
In fact the analog formula should correspond, again, to the 
logarithmic derivative of the zeta function as a product of L-functions:
$\sum_p \log L_p(\rho,\omega)$.

Since we are more and more tentative in statements, 
a further study from the point of view of theory of Hecke characters,
the role of the Teichmuler p-adic character and $\omega$-eigenvalues theory,
is required at this stage.

\section{Conclusions}
The primary goal of the presented investigations is to identify
an (pro-)algebraic structure underlying the Riemann spectrum,
dual to the POSet structure of primes:
$$(R\mbox{-}Spec=\{\theta_n\}_n\in N, \ ?\ ) \quad 
\leftrightarrow \quad (Z\mbox{-}Spec=\{p_n\}, <<).$$
In view of the Primes-Zeros duality \cite{Mazur-Primes}, 
it is natural to expect that this structure is of adelic origin.

Beyond the preliminary considerations presented above,
more specific questions need be addressed:

1) Is there a {\em primary spectrum} $\nu_p$ such that the 
{\em Riemann spectrum} \cite{Mazur-Primes} is derived from, 
like in the case of a multiplicative function?
Is this structure dual to the POSet structure of the primes?

2) Does the finite-infinite prime dichotomy apply beyond the
decomposition of the zeta function $\zeta\cdot \Gamma-factor=\xi(s)$?
We expect that the adelic duality of $Q$ implies an Algebraic Quantum Group duality
$$Q\to A_{alg} \to Hom(Q, T),$$
where $T=\sum Q_p/Z_p$ is the ``rational circle'' (roots of unity),
and that this reflects in the behavior of the zeta eigenvalues under
multiplication by $\log p$ discovered by Rademacher, Ford and Zaharescu
(the real scaling invaries the distribution).

We refer to this idea graphically as:
$$\xymatrix@R=1pc{
Primes: \ P \cup \infty & Zeros: \ \Theta + A \\
L-functions: \xi \quad 1/\Gamma & Adelic\ deccomposition: \theta_n=\sum_p\ \{\ \}_p + real
}$$
This adelic decomposition (of additive characters) 
could be related to the partial fraction decomposition of the rationals, 
which in turn is related to the Discrete Fourier Transform.

\vspace{.2in}
In conclusion, it seems that new ideas are needed for a 
break-through in the case of RH.
The general tendency is to move away from the analytic (NT) approach initiated 
by Riemann, which seems to have reached its limitations, 
and lean towards an adelic approach, notably the approach A. Connes is using \cite{Connes}.

A direct investigation of the Riemann zeros and of their imaginary parts 
is proposed, following \cite{Ford},
yet in connection with the POSet structure of prime numbers.
Since the zeros are dual to primes ($\delta_p$ basis,
either in the group algebra of the rationals as an algebraic quantum group,
or as distributions in the $\omega$-eigendistributions theory),
there should be a Fourier transform formula 
in terms of Hecke quasi-characters.

The work of Ford and Zaharescu is taken as a starting point,
with emphasis on the Rademacher case: the study of 
the change of the distribution of imaginary parts of zeta zeros
caused by the dilation operators $L_\alpha(\gamma_n)=\alpha\gamma_n$,
when $2\pi \alpha=\log p$.

The interpretation of the distribution function $g_\alpha$ as a 
local L-function is used to link it to the theory of $\omega$-eigendistributions
for such Hecke characters.
The approach makes use of the work of many other researchers providing various tools
to the research table: S. Kudla, R. Meyer, J. F. Burnol, A. Connes \cite{Kudla,Meyer,Burnol,Connes} 
etc. 

The theoretical understanding of Riemann Spectrum, in terms of 
p-adic characters, $\omega$-eigendistributions etc,
should be corroborated with a computer investigation:
1) determine several distributions $g_{\log p/2\pi}$ 
using SAGE or Mathematica, and look for correlations between $\lambda=\gamma_n/2\pi$
and primes $p$;
2) Investigate a p-adic interpretation of $\theta_n$, ``disguised''
under their decimal form, as a form of p-adic Haar wavelet analysis.

This proposed investigation was partially carried on, and reported in \cite{LI:StatsRSpec}.

Returning to the theoretical investigation,
since $Z\mbox{-}Spec=Spec(\C{A}b)$, 
and the POSet structure is induced by the loop ``functor'' $Aut_{Ab}$,
the adelic duality categorified invites to consider a Tannaka-Krein duality
between $\C{A}b_f$, the category of finite abelian groups 
and some representation category, probably involving the absolute Galois group,
the rationals as an agebraic quantum group, and crystal bases ...
a too big of a picture for the present author!

\section*{Acknowledgments}
I would like to thank G. Seelinger and P. Ding for their participation and help 
in our weekly seminar on the Riemann Hypothesis; and
for their patience regarding the illusive and intangible outcomes.


\end{document}